\documentclass[12pt]{amsart}
\usepackage{times}
\usepackage{amsfonts}
\usepackage{amsthm}
\usepackage{amsmath}
\advance \topmargin by -\headheight
\advance \topmargin by -\headsep
\evensidemargin \oddsidemargin
\newtheorem{theorem}{theorem}[section]
\newtheorem{lemma}[theorem]{Lemma}
\newtheorem*{hi}{Hyperfinite Inequality for $\delta_0$}
\newtheorem{corollary}[theorem]{Corollary}

\newtheorem{remark}[theorem]{Remark}

\newtheorem{definition}[theorem]{Definition}

\begin{document}

\title{A Hyperfinite inequality for free entropy dimension}

\author{Kenley Jung}

\address{Department of Mathematics, University of California,
Berkeley, CA 94720-3840,USA}

\dedicatory{For H-town}
\email{factor@math.berkeley.edu}
\thanks{Research supported by the NSF Graduate Fellowship Program}

\begin{abstract} If $X, Y, Z$ are finite sets of selfadjoint elements in
a tracial von Neumann algebra and $X$ generates a hyperfinite von
Neumann algebra, then $\delta_0(X \cup Y \cup Z) \leq \delta_0(X \cup Y) + 
\delta_0(X \cup Z)- \delta_0(X).$  We draw several corollaries from this 
inequality.   
\end{abstract}
\maketitle

In [11] Voiculescu describes the role of entropy in free probability.  
He discusses several problems in the area, one of which is the free 
entropy
dimension problem.  Free entropy dimension ([8], [9]) associates to 
an
$n$-tuple of selfadjoint operators, $X = \{x_1,\ldots, x_n\},$ in a
tracial von Neumann algebra $M$ a number $\delta_0(X)$ called the
(modified) free entropy dimension of $X.$ $\delta_0(X)$ is an asymptotic
Minkowski or packing dimension of sets of $n$-tuples of matrices which
model the behavior of $X.$ The free entropy dimension problem simply asks
whether $\delta_0(X) = \delta_0(Y)$ for any other $m$-tuple of selfadjoint
elements $Y$ satisfying $Y^{\prime \prime} = X^{\prime \prime}.$ It is
known from [10] that $\delta_0$ is an algebraic invariant, i.e.,
$\delta_0(X) = \delta_0(Y)$ if $X$ and $Y$ generate the same algebra.

The origin of this remark started with two
extremely special and highly tractable cases of this
problem, the first being: if $X, Y$ and $Z$  are finite sets of
selfadjoint elements in $M$ such that $X^{\prime \prime} = Z^{\prime
\prime}$ is hyperfinite, then is it true that 

\[ \delta_0(X \cup Y) = \delta_0(Y \cup Z)?\]

\noindent The second problem concerns invariance of $\delta_0$ over the
center: if $Y$ is an arbitrary set of selfadjoint elements in
$M$ and $y$ is any element in the center of $Y^{\prime \prime},$ then is
it true that

\[ \delta_0(Y \cup \{y\}) = \delta_0(Y)?\]

\noindent Both questions have affirmative answers and follow from a kind
of hyperfinite inequality for $\delta_0:$ If $X,Y, Z,$ are sets of
selfadjoint elements in $M$ and $X$ generates a hyperfinite von Neumann
algebra, then

\[ \delta_0(X \cup Y \cup Z) \leq \delta_0(X \cup Y) + \delta_0(X \cup Z) 
- \delta_0(X).
\]

\noindent Related inequalities of this nature can be found in
Gaboriau's work on the cost of equivalence relations [3].  The proof
of the microstates inequality above is an application of the work in
[5] paired with the packing formulation of $\delta_0$ in [6].

This remark has four sections.  The first is a short list of assumptions.
Motivated by the recent work of Belinschi and Bercovici ([1]), the second
presents a slightly simpler formulation of $\delta_0$ where the operator
cutoff constants are removed.  The third section presents the hyperfinite
inequality.  The fourth and last section presents several corollaries, two
of which answer the special invariance questions posed above.  It is also
shown that $\delta_0$ shares a certain property of the Connes-Shlyakhtenko
dimension, $\bigtriangleup.$ In [2] Connes and Shlyakhtenko show that if
$F$ is a finite set of selfadjoints in $M$ with $F^{\prime \prime}$ having
diffuse center, then $\bigtriangleup(F) = 1.$ We show that $\delta_0$
satisfies the same property provided that $F$ is assumed to have finite
dimensional approximants.

\section{Preliminaries}

Throughout suppose $M$ is a von Neumann algebra with a normal,
tracial state $\varphi.$ For any $n \in \mathbb N, |\cdot|_2$ denotes
the norm on $(M^{sa}_k(\mathbb C))^n$ given by $|(x_1,\ldots, x_n)|_2
= ( \sum_{j=1}^n tr_k(x_j^2))^{\frac{1}{2}}$ where $tr_k$ is the
tracial state on the $k\times k$ complex matrices, and
$|\cdot|_{\infty}$ denotes the operator norm.  $U_k$ denotes the $k
\times k$ unitary matrices.  For any $k,n \in \mathbb N, u \in U_k$
and $x = (x_1,\ldots, x_n) \in(M^{sa}_k(\mathbb C))^n,$ define $uxu^*
= (u x_1 u^*, \ldots, u x_n u^*).$ We will maintain the notation
introduced in [6], [8], and [9].  If $F = \{a_1, \ldots, a_n\}$ is a
finite set of selfadjoint elements in $M,$ we abbreviate
$\Gamma_R(a_1, \ldots, a_n;m,k,\gamma)$ by $\Gamma_R(F;m,k,\gamma)$
and in a similar way we write the associated microstate sets and
quantities introduced in [6], [8], and [9]: $\delta_0(F), \mathbb
P_{\epsilon}(F), \mathbb K_{\epsilon}(F).$ Also, if $G = \{b_1,
\ldots, b_p\}$ is another finite set of selfadjoint elements in $M,$
then we denote by $\Gamma_R(F \cup G;m,k, \gamma)$ the set
$\Gamma_R(a_1, \ldots, a_n, b_1, \ldots, b_p;m,k,\gamma)$ and write
all the associated microstate quantities $\delta_0(F \cup G), \mathbb
P_{\epsilon}(F \cup G),$ and $\mathbb K_{\epsilon}(F \cup G)$ with
respect to $\Gamma_R(F \cup G;m,k,\gamma).$ Finally,
$\Gamma(F;m,k,\gamma)$ will denote the set of all $k\times k$
microstates (no restrictions on the operator norms)  with degree of
approximation $(m,\gamma)$

\section{Cutoff Constants}

Recall that in [1] Belinschi and Bercovici have lifted the
operator norm cutoff constants in the definition of
$\chi.$  In other words, if $\chi(X)$ is the normal
definition as conceived of by Voiculescu, and
$\chi_{\infty}(X)$ is the quantity obtained by
replacing the microstate spaces
$\Gamma_R(X;m,k,\gamma)$ with $\Gamma(X;m,k,\gamma),$
then Belinschi and Bercovici showed that one always has

\[ \chi(X) = \chi_{\infty}(X).\]

We want to show the same thing for the packing formulation of
$\delta_0.$ Voiculescu defined $\delta_0(X)$ by

\[ \delta_0(X) = n + \limsup_{\epsilon \rightarrow 0} 
\frac{\chi(x_1 + \epsilon s_1, \ldots, x_n+ \epsilon   
s_n : s_1,\ldots, s_n)}{|\log \epsilon|} \]

\noindent where ${s_1,\ldots,s_n}$ is a semicircular 
family free with respect to $X.$ One can actually also
define $\delta_0(X)$ in terms of $\epsilon$ metric
packings.  One associates to each $\epsilon >0$ an   
asymptotic $\epsilon$ packing number $\mathbb 
P_{\epsilon}(X)$ and $\epsilon$ covering number
$\mathbb K_{\epsilon}(X).$ These definitions also make
use of cutoff constants.  Let's recall the definitions.
For any metric space $\Omega$ and $\epsilon >0$ denote
by $P_{\epsilon}(\Omega)$ the maximum number in a
collection of mutually disjoint open $\epsilon$ balls
of $\Omega$ and by $K_{\epsilon}(\Omega)$ the minimum 
number of open $\epsilon$ balls required to cover 
$\Omega.$ In what follows all spaces are endowed with
the $|\cdot |_2$ metric.  Define successively:
\[ \mathbb P_{\epsilon,r}(X;m,\gamma) = \limsup_{k
\rightarrow \infty} k^{-2} \cdot
\log(P_{\epsilon}(\Gamma_R(X;m,k,\gamma)),\]  

\[ \mathbb P_{\epsilon, r}(X) = \inf\{\mathbb
P_{\epsilon, R}(X;m,\gamma): m \in \mathbb N, \gamma
>0 \},\]

\[ \mathbb P_{\epsilon}(X) = \sup_{R>0} \mathbb
P_{\epsilon, R}(X).\]
\noindent Similarly we define $\mathbb K_{\epsilon}(X)$
by replacing all the $P_{\epsilon}$ above with
$K_{\epsilon}.$  It was shown in [5] that

\[ \delta_0(X) = \limsup_{\epsilon \rightarrow 0}
\frac{\mathbb P_{\epsilon}(X)}{|\log \epsilon|} =
\limsup_{\epsilon \rightarrow 0} \frac{\mathbb
K_{\epsilon}(X)}{|\log \epsilon|} \]

\noindent Now define successively:

\[ \mathbb P_{\epsilon, \infty}(X;m,\gamma) =
\limsup_{k \rightarrow \infty} k^{-2} \cdot \log
(\Gamma(X;m,k,\gamma),\]

\[ \mathbb P_{\epsilon, \infty}(X) = \inf \{\mathbb
P_{\epsilon, \infty}(X;m, \gamma): m \in \mathbb N,
\gamma >0\}. \]

\noindent Similarly we define $\mathbb K_{\epsilon,
\infty}.$  We want to show that the packing formulation for $\delta_0$ holds
when $\mathbb P_{\epsilon}(X)$ is replaced with
$\mathbb P_{\epsilon, \infty}(X)$ and $\mathbb
K_{\epsilon}(X)$ is replaced with $\mathbb K_{\epsilon,
\infty}(X).$  Suppose throughout this section, that $X$ is a finite 
set of selfadjoint elements in $M$ and that $R \geq 1$ is a 
constant greater than or equal to the maximum of the operator norms 
of the elements of $X.$  We need one easy lemma which is undoubtedly 
known, but which we will prove out for completeness:

\begin{lemma} For $m_0 \in \mathbb N,$ and $\epsilon, \gamma_0 > 0$
there exists an $m \in \mathbb N$ and $\gamma >0$ such that if $\xi =
(\xi_1,\ldots, \xi_n) \in \Gamma(X;m,k,\gamma),$ then $|\xi -
F_R(\xi)|_2 < \epsilon$ and $(F_R(\xi_{i_1}) \cdots F_R(\xi_{i_p})) \in
\Gamma_R(X;m_0,k,\gamma_0)$ where $F_r: \mathbb R \rightarrow [-r, r]$
is the monotone function equal to the identity on $(-r,r)$ and
$F_r(\xi) = (F_r(\xi_1),\ldots,F_r(\xi_n)).$
\end{lemma}

\begin{proof} Denote by $C$ the maximum over all 
numbers of the form $|F_{i_1}(\xi_1) \cdots
F_{i_p}(\xi_p)|_2+1$ where $1 \leq p \leq m$ and 
$i_1,\ldots, i_p \in \{R, \infty\}$ (this
constant $C$ is used to satisfy the second
condition).  Now choose $\gamma <
\frac{\epsilon}{100Cmn}.$
Choose $m \in \mathbb N$ so large that $\frac{R^m +
\gamma}{(R+ \gamma)^m} < \frac{\epsilon^2}{100n^2(R^4 +
\gamma)}.$ Suppose
$ \xi = (\xi_1,\ldots, \xi_n) \in
\Gamma(X;m_1,k,\gamma_1)$ and denote by
$\lambda_{i1},
\ldots, \lambda_{ik}$ the eigenvalues of $\xi_i$ with
multiplicity.  Set $B_i = \{ j \in \mathbb N: 1\leq j
\leq k, |\lambda_{ij}| \geq R+ \gamma \}.$  We have:

\[ \#B_i \cdot (R + \gamma)^m \leq \sum_{j \in B_i}
|\lambda_{ij}|^m \leq \sum_{j=1}^k |\lambda_j|^m \leq  
|Tr(\xi_i^m)| < k(R^m + \gamma).\]
\noindent Consequently $\frac{\#B_i}{k} \leq \frac{R^m
+ \gamma}{(R+ \gamma)^m}.$  By the Cauchy-Schwarz
inequality

\[|\xi_i - F_{R+\gamma}(\xi_i)|_2^2 \leq \frac{1}{k}
\cdot \sum_{j \in B_i} |\lambda_{ij}|^2 \leq |a^2|_2
\cdot \left (\frac{\#B_i}{k} \right )^{\frac{1}{2}} <
(R^4 + \gamma)^{\frac{1}{2}} \cdot
\sqrt{\frac{\epsilon^2}{100n^2(R^4+ \gamma)}} <
\frac{\epsilon}{10n}.\]

\noindent $|F_{R + \gamma}(\xi_i) - F_R(\xi)|_2 \leq
\gamma$ whence it follows that $|\xi -F_R(\xi)|_2 <
\epsilon.$   To see that the second claim is satisfied
observe that for any $1 \leq p \leq m$ and $1 \leq i_1, \ldots, i_p
\leq n,$ Cauchy-Schwarz again yields

\[|tr_k(\xi_{i_1} \cdots \xi_{i_p}) - tr_k(F(\xi_{i_1}) \cdots F(\xi_{i_p}))| \leq C \cdot p \cdot |\xi -
F_r(\xi)|_2. \]

\noindent Because the $|\cdot|_2$ quantity on the  
right hand side can be made smaller than any given
$\epsilon > 0$ and $C$ and $p$ are constants bounded
from the get-go, we're done.
\end{proof}

\begin{lemma}

\begin{eqnarray*} \delta_0(X) = \limsup_{\epsilon
\rightarrow 0} \frac{\mathbb P_{\epsilon,
\infty}(X)}{|\log \epsilon|} & = &\limsup_{\epsilon
\rightarrow 0} \frac{\mathbb K_{\epsilon,
\infty}(X)}{|\log \epsilon|} \\ & = & \limsup_{\epsilon
\rightarrow 0} \frac{\mathbb P_{\epsilon}(X)}{|\log
\epsilon|} \\ & = & \limsup_{\epsilon \rightarrow
0} \frac{K_{\epsilon}(X)}{|\log \epsilon|}.
\end{eqnarray*}

\end{lemma}

\begin{proof} Suppose $m_0 \in \mathbb N$ and $\gamma_0
>0.$  There exist by Lemma 2.1 an $m \in \mathbb N$
and $\gamma >0$ such that if $\xi \in
\Gamma(X;m,k,\gamma),$ then $|\xi - F_R(\xi)|_2 <  
\epsilon$ and it can be also arranged so that
$F_R(\xi) \in \Gamma_R(X;m_0,k,\gamma_0).$  It follows
from this that $K_{2 \epsilon}(\Gamma(X;m,k,\gamma))
\leq K_{\epsilon}(\Gamma_R(X;m_0,k,\gamma_0)$ whence 

\[ \mathbb K_{2\epsilon, \infty}(X) \leq \mathbb K_{2  
\epsilon}
(X;m, \gamma) \leq
\mathbb K_{\epsilon, R}(X, m_0, \gamma_0).
\]

\noindent So $\mathbb K_{2\epsilon,
\infty}(X) \leq
\mathbb K_{\epsilon, R}(X) \leq \mathbb
K_{\epsilon}(X).$  Now clearly

\[ \delta_0(X) = \limsup_{\epsilon \rightarrow 0}
\frac{\mathbb K_{\epsilon}(X)}{|\log \epsilon|} \leq
\limsup_{\epsilon \rightarrow 0} \frac{\mathbb
K_{2 \epsilon, \infty}(X)}{|\log 2 \epsilon|} \leq
\limsup_{\epsilon \rightarrow 0} \frac{\mathbb
K_{\epsilon}(X)}{|\log 2 \epsilon|} = \delta_0(X). \]

\noindent $\mathbb P_{\epsilon, \infty}(X) \geq \mathbb
K_{2 \epsilon, \infty}(X) \geq \mathbb P_{4 \epsilon,
\infty}(X)$ so this completes the proof.
\end{proof}

\section{The Inequality}

Throughout assume $X, Y, Z$ and $F$ are finite sets of selfadjoint
elements in $M.$ Assume further that $X$ generates a hyperfinite von
Neumann algebra, an assumption we will restate for emphasis in some
of the corollaries.

\begin{definition}  Suppose for each $m \in \mathbb N$ and $\gamma >0,
\langle \xi_k \rangle_{k=1}^{\infty}$ is a sequence such that for large
enough $k, \xi_k \in \Gamma(X;m,k,\gamma).$  The set of 
microstates $\Xi(F;m,k,\gamma)$ for $F$ relative to the $\xi_k$ is

\[ \Xi(F;m,k,\gamma) = \{\eta : (\xi_k, \eta) \in \Gamma(X \cup F;m,
k,\gamma) \}.\]

\noindent Define successively for $\epsilon >0,$

\[ \mathbb K_{\epsilon}(\Xi(F;m,\gamma)) = \limsup_{k \rightarrow
\infty} k^{-2} \cdot \log K_{\epsilon}(\Xi(F;m,k,\gamma)),\]

\[ \mathbb K_{\epsilon}(\Xi(F)) = \inf \{ \mathbb K_{\epsilon}(F;m,
\gamma) : m \in \mathbb N, \gamma >0\}. \]

\noindent where the packing quantities are taken with respect to
$|\cdot |_2.$ In a similar fashion, we define $\mathbb
P_{\epsilon}(\Xi(F))$
by replacing the $K_{\epsilon}$ above with $P_{\epsilon}.$  

\end{definition}

\begin{lemma} Suppose for each $m \in \mathbb N$ and $\gamma >0$ we have a
sequence $\langle \xi_k \rangle_{k=1}^{\infty}$ satisfying $\xi_k \in
\Gamma(X;m,k,\gamma)$ for sufficiently large $k.$ Then,

\[ \delta_0(X \cup F) = \delta_0(X) + \limsup_{\epsilon \rightarrow 0}
\frac{\mathbb K_{\epsilon}(\Xi(F))}{|\log \epsilon|}. \]

\end{lemma}

\begin{proof} First we show that the left hand side is greater than
or equal to the right hand side.  Suppose $t >0$ is given. By [5] and
[6] there exists an $\epsilon_0 >0$ such that for all $\epsilon_0 >
\epsilon >0$ and any $m \in \mathbb N, \gamma >0, \liminf_{k
\rightarrow \infty} P_{\epsilon} (\Gamma(X;m,k,\gamma)) >
(\delta_0(X) - t) |\log 2 \epsilon|.$ Now suppose $m \in \mathbb N$
and $\gamma >0$ are fixed.  Consider the $\langle \xi_k
\rangle_{k=1}^{\infty}$ corresponding to the $m$ and $\gamma.$
Because the von Neumann algebra generated by $X$ is hyperfinite, by
[5] I can find a set of unitaries $\langle v_{\lambda k}
\rangle_{\lambda \in \Lambda_k}$ such that $\langle v_{\lambda k}
\xi_k v_{\lambda k}^* \rangle_{\lambda \in \Lambda_k}$ is an
$\epsilon$ separated set with respect to the $| \cdot |_2$ norm and
$\liminf_{k\rightarrow} k^{-2} \cdot \log \#\Lambda_k > (\delta_0(X)
- t) | \log 2 \epsilon|.$ For each $k$ pick an $\epsilon$ separated
subset $\langle \zeta_k \rangle_{k=1}^{\infty}$ of minimal
cardinality for $\Xi(F;m,k,\gamma)$ (the set of microstates for $F$
relative to $\xi_k$).  Now it is manifest that

\[ \langle (v_{\lambda k} \xi_k v_{\lambda k}^*, v_{\lambda k} \zeta_{jk}
v_{\lambda k}^*) \rangle_{(\lambda,j) \in \Lambda_k \times J_k} \]

\noindent is a subset of $\Gamma(X \cup F;m,k,\gamma)$ and moreover, it 
is easily checked that this set is $\epsilon$-separated with respect to 
the $|\cdot |_2$ norm.  Hence,

\begin{eqnarray*} \mathbb P_{\epsilon, \infty} (X \cup F;m, \gamma) &
\geq & \limsup_{k\rightarrow \infty} k^{-2} \cdot \log (\# \Lambda_k
\cdot \#J_k)\\ &\geq & \liminf_{k \rightarrow \infty} k^{-2} \cdot
\log \#\Lambda_k + \limsup_{k \rightarrow \infty} k^{-2} \cdot \log
P_{\epsilon} (\Xi(F;m,k, \gamma) \\ & \geq & (\delta_0(X) - t) | \log
2 \epsilon| + \limsup_{k \rightarrow \infty} k^{-2} \cdot \log
P_{\epsilon}(\Xi(F; m,k, \gamma) \end{eqnarray*}

\noindent so that for $\epsilon_0 > \epsilon >0$

\[ \mathbb P_{\epsilon, \infty}(X \cup F) \geq  (\delta_0(X) - t) 
|\log 2 \epsilon| + \mathbb P_{\epsilon}(\Xi(F)) \]

\noindent Using the packing formulation of $\delta_0$ in [6], the 
fact that for any metric space $\Omega, P_{\epsilon}(\Omega) \geq 
K_{2
\epsilon}(\Omega) \geq P_{4 \epsilon}(\Omega),$ and dividing by $|\log
\epsilon|$ and taking $\limsup_{\epsilon \rightarrow 0}$ on both sides
gives

\begin{eqnarray*} \delta_0(X \cup F) = \limsup_{\epsilon \rightarrow 0}
\frac{\mathbb P_{\epsilon, \infty}(X \cup F)}{|\log \epsilon|} &\geq&
\limsup_{\epsilon \rightarrow0} \frac{ (\delta_0(X) - t) |\log 2 \epsilon|
+ \mathbb P_{\epsilon}(\Xi(F))}{|\log \epsilon|} \\ & = & \delta_0(X) -
t +\limsup_{\epsilon \rightarrow 0} \frac{\mathbb
P_{\epsilon}(\Xi(F))}{|\log \epsilon|} \\ & = & \delta_0(X) - t
+\limsup_{\epsilon \rightarrow 0} \frac{\mathbb
K_{\epsilon}(\Xi(F))}{|\log \epsilon|}. \\ \end{eqnarray*}

\noindent $t >0$ being arbitrary we have 

\[ \delta_0(X \cup F) \geq \delta_0(X) + \limsup_{\epsilon \rightarrow 0}
\frac{\mathbb P_{\epsilon}(\Xi(F))}{|\log \epsilon|}.\]

For the reverse inequality by [6] there are $C, \epsilon_0 >0$ such
that for $\epsilon_0 > \epsilon >0$ and for any $k \in \mathbb N$ and
tractable subgroup $H$ of $U_k$ (in the sense of [4]) there exists an
$\epsilon$-net for $U_k/H$ with respect to the quotient metric
induced by $| \cdot|_{\infty}$ with cardinality no greater than
$\left(\frac{C}{\epsilon}\right)^{\dim(U_k/H)}.$ Write $R$ for the
maximum of the operator norms of the elements in $X.$ Suppose $m \in
\mathbb N$ and $\gamma >0.$ Observe that there exists $\epsilon > r
>0$ so small that if $(\xi, \eta) \in \Gamma_R(X \cup
F;m,k,\gamma/2)$ and $|(\xi, \eta) - (x, a)|_2 < r,$ then $(x, a) \in
\Gamma(X \cup F;m,k,\gamma).$ There also exist $m_1 \in \mathbb N$
and $\gamma_1 >0$ such that if $\xi, x \in
\Gamma(X;m_1,k,\gamma_1),$ then there exists a $u \in U_k$
satisfying $|u\xi u^* -x|_2 < r.$ Set $m_2 =m + m_1$ and $\gamma_2 =
\min \{ \gamma/2, \gamma_1).$

By [6] I can find a sequence $\langle \xi_k \rangle_{k=1}^{\infty}$
such that for sufficiently large $k$ $\xi_k \in
\Gamma_R(X;m_1,k,\gamma_1)$ and $\dim \xi_k^{\prime} \geq k^2 (1
-\delta_0(X)).$ Consider the associated $\Xi(F;m,k,\gamma)$ and for
each $k$ find an $\epsilon$-net $\langle \eta_{jk} \rangle_{j \in
J_k}$ for $\Xi(F;m,k,\gamma)$ with respect to $| \cdot |_2$ of
minimum cardinality.  Define $H_k$ to be the unitary group of
$\xi_k^{\prime}.$ I can find for each $k$ large enough a set of
unitaries $\langle u_{gk} \rangle_{g \in G_k}$ such that their images
in $U_k/H_k$ is an $\epsilon$-net with respect to the quotient metric
induced by $| \cdot |_{\infty}$ and such that

\[ \# \Lambda_k \leq \left (\frac{C}{\epsilon} \right)^{\delta_0(X)
k^2}.\]

\noindent Consider 

\[ \langle (u_{gk} \xi_k u_{gk}^*, u_{gk}\eta_{jk}u_{gk}^*) 
\rangle_{(g,j) \in G_k \times J_k}.\]

\noindent I claim that this set is a $5 \epsilon R$-net for $\Gamma(X 
\cup F;m_2,k,\gamma_2).$

To see this suppose $(\xi, \eta) \in \Gamma(X \cup F;m_2, k, \gamma_2).$
By the selection of $m_1$ and $\gamma_1$ there exists a $u \in U_k$ such
that $|u^* \xi_k u - \xi|_2 < r_.$ Taking into account the stipulation on
$r$ this implies that ($u^* \xi_k u, \eta) \in \Gamma(X ;m,k,\gamma)
\Longleftrightarrow (\xi_k, u \eta u^*) \in \Gamma(X \cup
F;m,k,\gamma),$ whence $u \eta u^* \in \Xi(F;m,k,\gamma).$ There exists
an $g \in G_k$ and an $h \in H_k$ such that $|u - u_{gk} h|_{\infty} <
\epsilon.$ Consequently,

\[ |u_{gk}\xi_k u_{gk}^* - \xi|_2 = |u_{gk} h \xi_k h^* u_{gk} - \xi|_2
\leq 2\epsilon R + |u \xi_k u^* - \xi|_2 \leq 3 \epsilon R. \]

\noindent Now $u \eta u^* \in \Xi(F;m,k,\gamma)$ so there exists a $j \in
J_k$ such that $|\eta_{jk} - u \eta u^* |_2 < \epsilon.$ Because
$\Xi(F;m,k,\gamma)$ is invariant under the action of $H_k$ it follows that
$h \eta_{jk} h^* \in \Xi(F;m,k,\gamma)$ and so there exists an $\ell \in
J_k$ such that $|\eta_{\ell k} - h \eta_{jk} h^*|_2 < \epsilon.$ So again
we have

\[ |u_{gk} \eta_{\ell k} u_{gk}^* - \eta|_2 < |u_{gk} h \eta_{jk} h^*
u_{gk}^* - \eta|_2 + \epsilon < |u \eta_{jk} u^* - \eta|_2 + 3 \epsilon R 
< 4 \epsilon R.
\]

\noindent and we have the desired claim.

It follows that

\begin{eqnarray*} \mathbb K_{5\epsilon R, \infty}(X \cup F;m_2,
\gamma_2) & \leq & \limsup_{k \rightarrow \infty} k^{-2} \cdot \log(
\# G_k \cdot \# J_k ) \\ & \leq & \log C + \delta_0(X) \cdot |\log
\epsilon| + \limsup_{k \rightarrow \infty} k^{-2} \cdot \log
K_{\epsilon}(\Xi(F;m,k,\gamma)).  \end{eqnarray*}

\noindent Given any $m \in \mathbb N$ and $\gamma >0$ we produced $m_2 \in
\mathbb N$ and $\gamma_2 >0$ so that the above inequality holds for 
$0 < \epsilon < \epsilon_0.$  Thus  

\[ \mathbb K_{5 \epsilon, \infty} (X \cup F) \leq \log C + 
\delta_0(X) \cdot
|\log\epsilon| + \mathbb K_{\epsilon}(\Xi(F))). \]

\noindent Taking $\limsup_{\epsilon \rightarrow 0}$ on both sides and
again using the packing formulation of $\delta_0$ in [6] as well as 
Lemma 2.2. we have 

\begin{eqnarray*} \delta_0(X \cup F)  = \limsup_{\epsilon \rightarrow 0}
\frac{\mathbb K_{5 \epsilon R, \infty}(X \cup  F)}{|\log 5 \epsilon 
R|} & \leq 
& \limsup_{\epsilon \rightarrow 0} \frac{\log C + \delta_0(X) \cdot |\log
\epsilon| + \mathbb K_{\epsilon}(\Xi(F))}{|\log 5 \epsilon R|} \\ & = &
\delta_0(X) + \limsup_{\epsilon \rightarrow 0} \frac{\mathbb K_{\epsilon}(\Xi(F))}{|\log \epsilon|} \\ \end{eqnarray*}. \end{proof}

\begin{hi} If $X^{\prime \prime}$ is hyperfinite, then

\[\delta_0(X \cup Y \cup Z) \leq \delta_0(X \cup Y) + \delta_0(X \cup Z) - 
\delta_0(X).\]
\end{hi}

\begin{proof} $X$ has finite dimensional approximants, so for each $m \in
\mathbb N$ and $\gamma >0$ we can find sequences $\langle \xi
\rangle_{k=1}^{\infty}$ satisfying the conditions of Lemma 2.2 and
consider all relative microstates with respect to these fixed sequences.  
For each $k, \Xi(Y \cup Z;m,k,\gamma) \subset \Xi(Y;m,k,\gamma) 
\times
\Xi(Z;m,k,\gamma)$ so that
 
\[K_{2 \epsilon}(\Xi(Y \cup Z;m,k,\gamma)) \leq
K_{\epsilon}(\Xi(Y;m,k,\gamma)) \cdot
K_{\epsilon}(\Xi(Z;m,k,\gamma)).\]

\noindent It follows that 

\[ \limsup_{\epsilon \rightarrow 0} \frac{\mathbb K_{\epsilon}(\Xi(Y
\cup Z))}{|\log \epsilon|} \leq \limsup_{\epsilon \rightarrow 0}
\frac{\mathbb K_{\epsilon}(\Xi(Y))}{|\log \epsilon|} +
\limsup_{\epsilon\rightarrow 0} \frac{\mathbb
K_{\epsilon}(\Xi(Z))}{|\log \epsilon|}.\]

\noindent By the preceding lemma and the inequality above

\begin{eqnarray*} \delta_0(X \cup Y \cup Z) & = & \delta_0(X) +
\limsup_{\epsilon\rightarrow 0} \frac{\mathbb
K_{\epsilon}(\Xi( Y \cup
 Z))}{|\log \epsilon|} \\ & \leq & \delta_0(X) + \limsup_{\epsilon \rightarrow 0}
\frac{\mathbb K_{\epsilon}(\Xi(Y))}{|\log \epsilon|} +\limsup_{\epsilon
\rightarrow 0} \frac{\mathbb K_{\epsilon}(\Xi(Z))}{|\log \epsilon|} \\ & = &
\delta_0(X \cup Y) + \delta_0(X \cup Z) -\delta_0(X). \end{eqnarray*} \end{proof}

\begin{remark} The hyperfinite assumption on $X^{\prime \prime}$ is
necessary.  To see this consider the group inclusion $\mathbb F_3 \subset
\mathbb F_2 \subset \mathbb F_3$ where $\mathbb F_n$ is the free group on
$n$ generators.  On the von Neumann algebra level this translates to
$L(\mathbb F_3) \simeq M_1 \subset L(\mathbb F_2) \simeq M_2 \subset
L(\mathbb F_3) \simeq M_3.$ Take $X, Y$ and $Z$ to be the canonical sets
of freely independent semicirculars associated to $M_1, M_2,$ and $M_3,$
respectively.  Then $\delta_0(X \cup Y) + \delta_0(X \cup Z) - \delta_0(X)
= 3 + 2 - 3 = 2 < 3 =\delta_0(X \cup Y \cup Z).$ \end{remark}

\section{Five Corollaries}

In this section $X, Y,$ and $Z$ are again finite sets of selfadjoint
elements in $M.$  Here are some corollaries of the hyperfinite inequality
for $\delta_0$:

\begin{corollary} Suppose $X^{\prime \prime}$ is hyperfinite.  Assume one
of the following holds:

\begin{itemize}
\item $Z \subset X^{\prime \prime}.$
\item $X^{\prime \prime}$ is diffuse, $ \delta_0(X\cup Z) \leq 1,$ and $Z
\subset (X \cup Y)^{\prime \prime}.$
\end{itemize}

\noindent Then $\delta_0(X \cup Y) = \delta_0(X \cup Y \cup Z).$
\end{corollary}
		
\begin{proof} In either of the two cases $Z$ is contained in the von
Neumann algebra generated by $X$ and $Y$ so by [10] $\delta_0(X \cup 
Y) 
\leq \delta_0(X \cup Y \cup Z).$  For the reverse inequality observe that 
either situations imply $\delta_0(X \cup Z) = \delta_0(X).$  This follows 
in the first case from invariance of $\delta_0$ for hyperfinite von 
Neumann algebras ([5]).  In the second case we have by assumption and 
hyperfinite
monotonicity that $1 \geq \delta_0(X \cup Z) \geq \delta_0(X) \geq 1.$
In either cases $\delta_0(X \cup Z) = \delta_0(X)$ so by the hyperfinite
inequality,

\[ \delta_0(X \cup Y \cup Z) \leq \delta_0(X \cup Y) + \delta_0(X \cup Z) 
- \delta_0(X) = \delta_0(X \cup Y).\] 

\noindent Thus, $\delta_0(X \cup Y) = \delta_0(X \cup Y \cup Z).$
\end{proof}

\begin{corollary} If $X^{\prime \prime} = Z^{\prime \prime}$
is hyperfinite, then $\delta_0(X \cup Y) = \delta_0(Y \cup Z).$
\end{corollary}

\begin{corollary} If $y= y^*$ lies in the center of the von Neumann
algebra generated by $Y,$ 
then

\[ \delta_0(Y \cup \{y\}) = \delta_0(Y).\]
\end{corollary}

\begin{proof} Again by [10] $\delta_0(Y) \leq \delta_0(Y \cup \{y\} 
).$ For 
the reverse inequality set $\alpha = \sup \{ \delta_0(x): x=x^* \in 
Y^{\prime \prime}\}$ (actually the supremum is achieved but I won't need 
that).  
Suppose $\epsilon >0.$ Find $x = x^* \in Y^{\prime \prime}$ such that
$\alpha - \epsilon < \delta_0(x).$ Take a sequence $\langle x_k
\rangle_{k=1}^{\infty}$ such that for each $k,$ $x_k = x_k^*$ lies in the
$*$-algebra generated by $Y$ and such that $x_k \rightarrow x$ strongly.  
Now for every $k$ there exists an $a_k = a_k^*$ such that the von Neumann
algebra generated by $a_k$ is equal to the von Neumann algebra generated
by $x_k$ and $y$ and thus $\delta_0(a_k) = \delta_0(x_k,y).$ Using the 
fact that $\delta_0$ is an algebraic invariant we have by the hyperfinite
inequality for $\delta_0$ 

\begin{eqnarray*} \delta_0(Y \cup \{y\}) = \delta_0(\{x_k\}\cup Y\cup
\{y\}) & \leq & \delta_0(\{x_k\} \cup Y)+ \delta_0(x_k,y) - \delta_0(x_k)  
\\ & = & \delta_0(Y) + \delta_0(a_k) - \delta_0(x_k) \\ & \leq &
\delta_0(Y) + \alpha - \delta_0(x_k). \\ \end{eqnarray*}

\noindent Forcing $k \rightarrow \infty$ and using the fact that
$\liminf_{k \rightarrow \infty} \delta_0(x_k) \geq \delta_0(x)$ (by 
[8])
we have that $\delta_0(Y \cup \{y\}) \leq \delta_0(Y) + \alpha - (\alpha -
\epsilon) = \delta_0(Y) + \epsilon.$  $\epsilon >0$ being arbitrary,
$\delta_0(Y\cup \{y\}) \leq \delta_0(Y).$  Thus, $\delta_0(Y \cup \{y\}) = 
\delta_0(Y).$
\end{proof}

\begin{corollary} Suppose $x = x^* \in M, \delta_0(x,Y) = \alpha,
\delta_0(Z) = \beta, \{x\} \cup Y \subset Z^{\prime \prime},$ and $Z =
\{z_1,\ldots, z_n\}.$  Then

\[ \beta - \alpha + n \cdot \delta_0(x) \leq \sum_{j=1}^n \delta_0(x,
z_j).\]

\noindent Thus if $n < \beta - \alpha + n \cdot \delta_0(x),$ then for
some $j, 1 < \delta_0(x, y_j).$  In particular, if $Z$ consists of $2 \leq
\beta \in \mathbb N$ freely independent semicircular elements, $Z = \{s_1,
\ldots, s_{\beta}\}$ and $x$ is any self-adjoint 
element in $Z^{\prime \prime}$ with no atoms, then for some $1 \leq j \leq
\beta, 1 < \delta_0(x, s_j).$
\end{corollary}

\begin{proof} $\{x\} \cup Y \subset Z^{\prime \prime}$ so by [9] and 
the hyperfinite inequality

\[ \beta = \delta_0(Z) \leq \delta_0(\{x\} \cup Y \cup \{z_1,\ldots, 
z_n\}) \leq \delta_0(\{x\} \cup
Y \cup \{z_1, \ldots, z_{n-1}\} ) + \delta_0(x,z_n) - \delta_0(x).\]

\noindent Repeating this $n$ times we arrive at

\[ \beta \leq \delta_0(\{x\} \cup Y) + \sum_{j=1}^n \delta_0(x, z_j) - n 
\cdot \delta_0(x) = \alpha - n \cdot \delta_0(x) + \sum_{j=1}^n
\delta_0(x,z_j),\]

\noindent whence $\beta - \alpha + n \cdot \delta_0(x) \leq \sum_{j=1}^n
\delta_0(x, z_j).$  Everything else is obvious.
\end{proof}

\begin{remark} Recall from [4] that for a finite set of selfadjoint
elements $F$ in $M,$ if $\delta_0(F) >1,$ then the von Neumann algebra
generated by $F$ cannot be generated by a sequence of Haar unitaries
$\langle u_j \rangle_{j=1}^{s}$ satisfying the condition $u_{j+1} u_j
u_{j+1}^* \in \{ u_1, \ldots, u_j\}^{\prime \prime}.$ In particular
$F^{\prime \prime}$ is prime and has no Cartan subalgebras.  Thus, in
the context of Corollary 4.4 for some $j, \{x, s_j\}^{\prime \prime}$
is prime and has no Cartan subalgebra.  \end{remark}

We conclude with a microstates analogue of a property of the 
Connes-Shlyakhtenko dimension $\bigtriangleup :$

\begin{corollary} If $X = \{x, x_1,\ldots, x_n\},$ $x x_i = x_i x$ for all 
$1 \leq i \leq n,$ and the spectrum of $x$ is diffuse, then $\delta_0(X) 
\leq 1.$  If in addition, $X$ has finite dimensional approximants, then 
$\delta_0(X) =1.$  Consequently, if $F$ is a finite set of selfadjoint 
elements in $M$ which has finite dimensional approximants and such that 
the von Neumann algebra generated by $F$ has diffuse center, then 
$\delta_0(F) =1.$  \end{corollary}

\begin{proof} By the hyperfinite inequality for $\delta_0,$ the 
diffuseness of $x,$ and [9]

\[ \delta_0(X) \leq \delta_0(x, \ldots, x_{n-1}) + \delta_0(x, x_n) - 
\delta_0(x) \leq \delta_0(x, \ldots, x_{n-1}) + 1 -1 = \delta_0(x, 
x_1,\ldots, x_{n-1}).\]

\noindent Continuing inductively we have $\delta_0(X) \leq \delta_0(x) =1$ 
as promised.  If $X$ has finite dimensional approximants, then by [5] 
$\delta_0(X) \geq \delta_0(x) =1$ and consequently, $\delta_0(X) =1.$  The 
claim concerning $F$ is immediate.

\end{proof}

\noindent{\it Acknowledgements.} A large part of this research was
conducted at UCLA and I thank Dimitri Shlyakhtenko for reading a
preliminary version of this remark, pointing out the relation between the
hyperfinite inequality with the work of Gaboriau, proposing the right
argument for the cutoff constant comment, and suggesting Corollary 4.6.  I
also thank him and the UCLA mathematics department for their warm
hospitality.

\end{document}